\documentclass[AMA,Times1COL]{WileyNJDv5} 

\usepackage{amsmath}
\usepackage{amsthm}
\usepackage[short]{optidef}
\usepackage{tikz}
\usepackage{orcidlink}

\usepackage{custom_SCL}

\articletype{Research Article}%

\received{Date Month Year}
\revised{Date Month Year}
\accepted{Date Month Year}
\journal{Journal}
\volume{00}
\copyyear{2024}
\startpage{1}

\begin{document}

\title{A Convex Optimization Approach to Compute Trapping Regions for Lossless Quadratic Systems}
\titlemark{A Convex Optimization Approach to Compute Trapping Regions for Lossless Quadratic Systems}

\author[1]{Shih-Chi Liao}

\author[2]{A. Leonid Heide}

\author[2]{Maziar S. Hemati}

\author[1]{Peter J. Seiler}

\authormark{LIAO, HEIDE, HEMATI, and SEILER}

\address[1]{\orgdiv{Electircal and Computer Engineering}, \orgname{University of Michigan, Ann Arbor}, \orgaddress{\state{Michigan}, \country{USA}}}

\address[2]{\orgdiv{Aerospace Engineering and Mechanics}, \orgname{University of Minnesota}, \orgaddress{\state{Minnesota}, \country{USA}}}

\corres{Shih-Chi Liao, Electrical and Computer Engineering, University of Michigan, Ann Arbor.  \email{shihchil@umich.edu}}


\fundingInfo{This material is based upon work supported by the Army Research Office under grant number W911NF-20-1-0156 and the Air Force Office of Scientific Research under grant number FA9550-21-1-0434. }

\abstract[Abstract]{%
Quadratic systems with lossless quadratic terms arise in many applications, including models of atmosphere and incompressible fluid flows. Such systems have a trapping region if all trajectories eventually converge to and stay within a bounded set. Conditions for the existence and characterization of trapping regions have been established in prior works for boundedness analysis. However, prior solutions have used non-convex optimization methods, resulting in conservative estimates. 
In this paper, we build on this prior work and provide a convex semidefinite programming condition for the existence of a trapping region. The condition allows precise verification or falsification of the existence of a trapping region. If a trapping region exists, then we provide a second semidefinite program to compute the least conservative trapping region in the form of a ball. Two low-dimensional systems are provided as examples to illustrate the results. A third high-dimensional example is also included to demonstrate that the computation required for the analysis can be scaled to systems of up to $\sim O(100)$ states. 
The proposed method provides a precise and computationally efficient numerical approach for computing trapping regions. We anticipate this work will benefit future studies on modeling and control of lossless quadratic dynamical systems. 
}


\jnlcitation{
\cname{%
    \author{Liao S},
    \author{Heide}, 
    \author{Hemati}, and
    \author{Seiler P}}.
\ctitle{A Convex Optimization Approach to Compute Trapping Regions for Lossless Quadratic Systems} 
\cjournal{\it J International Journal of Robust and Nonlinear Control} 
\cvol{....}.
}

\maketitle



\section{Introduction}


Dynamical systems with lossless quadratic nonlinearities arise in models of many physical phenomena, including those of the atmosphere~\cite{lorenz1963deterministic} and incompressible fluid flows~\cite{schlegel2015long}. One important feature of such systems is the long-term boundedness, i.e., the energy of the system remains finite. As most of the concerned systems are inherently bounded by their nature, models for such systems should also exhibit the same boundedness. Similarly, closed-loop control of lossless systems has been studied in the literature~\cite{MUSHTAQ2024ConvexSOF}. In such cases, the boundedness can be used to maintain the closed-loop system in a certain operation condition~\cite{khalil2002nonlinear}. Hence, methods to verify boundedness are required for this class of models. One straightforward approach is establishing an equilibrium's global stability, entailing all trajectories to converge and be bounded. However, a system can have an unstable equilibrium yet still be bounded, e.g., the Lorenz system~\cite{lorenz1963deterministic}. Thus, the typical Lyapunov stability analysis is often too strict for such systems, and a dedicated boundedness analysis is required.

A trapping region is a bounded set that all trajectories eventually converge into and stay within. This notion has been proposed in the literature as an approach to analyze boundedness. Lorenz~\cite{lorenz1963deterministic} proposed an energy function as a certificate for the existence of a trapping region to study the boundedness of the chaotic three-state system. Schlegel and Noack~\cite{schlegel2015long} formalized this idea with a necessary and sufficient condition for the existence of a spherical trapping region. They applied the condition to demonstrate the boundedness of several reduced-order models of incompressible fluid flows with $\sim O(10)$ state variables. Furthermore, these trapping region conditions have been exploited in physics-informed system identification algorithms~\cite{Kaptanoglu2021, goyal2023guaranteed} to obtain long-term bounded models of fluid flows with $\sim O(10)$ state variables from data. 

While trapping regions provide a certificate of boundedness, verifying the existence of a trapping region is non-trivial. Several approaches are discussed in Section 7 of Schlegel and Noack~\cite{schlegel2015long}. However, they are either limited by the state dimension of systems (as they involve an algebraic condition on the eigenvalues of the system matrices) or involve non-convex optimization relying on a good initial seed for the algorithm. The approaches offer no guarantees despite using physical understanding and simulations to aid the analysis. Thus the approaches may fail to find a trapping region even if one exists. Furthermore, even if the existence of a trapping region is verified, Schlegel and Noack provide a conservative radius estimation of the spherical trapping region. These limitations restrict the application of the boundedness analysis using trapping regions. 

This paper provides three contributions to the trapping region analysis. First, a convex semidefinite programming (SDP) condition is proposed as the necessary and sufficient condition for the system to have a trapping region. This convex condition provides theoretical guarantees and corresponding efficient numerical solutions. Hence, readily available SDP solvers can be used to verify or falsify the existence of a trapping region. Second, a quadratically constrained quadratic program (QCQP) is derived to compute the least conservative radius of the spherical trapping region if one exists. This QCQP can be solved exactly by its SDP dual problem. Also, the solution of the dual SDP characterizes the set of points with instantaneously non-decreasing energy on the boundary of the trapping region. This set of points can be useful for further study via simulations. Lastly, the proposed SDP-based trapping region analysis is demonstrated by three examples. The first two examples are low-dimensional systems to illustrate the concepts and the analysis process. The last example highlights that the proposed method can scale to $\sim O(100)$ states. The SDP-based method can easily verify the existence of the trapping region and provide the least conservative estimate for the trapping region. This convex formulation can potentially benefit the relevant modeling~\cite{Kaptanoglu2021, goyal2023guaranteed} and control~\cite{MUSHTAQ2024ConvexSOF} approaches for lossless quadratic systems.


The rest of the paper is structured in four sections. Section~\ref{sect:Prel} formulates the problem of this study and presents relevant prior work. Section~\ref{sect:Verif} presents the proposed SDP-based trapping region analysis, including verification, computing the least conservative radius, and characterizing critical points. Section~\ref{sect:Examp} provides numerical examples to illustrate the process and computational load of the analysis. Lastly, Section~\ref{sect:Conclusion} summarizes the paper and discusses possible future directions.

\section{Preliminaries}     \label{sect:Prel}



\subsection{Lossless Quadratic Systems} \label{sect:Prel_EnergyPresQuadSys}

Let $c\in\Rn{n}$, $L \in\Rn{n\times n}$, and symmetric matrices $\{Q^{(i)}\}_{i=1}^n \subset \Rn{n\times n}$ be given. Let $f: \Rn{n} \to \Rn{n}$ be a quadratic function with $i^{th}$ entry defined by $f_i(x) = x^\top Q^{(i)} x$ for $i=1\ldots n$. Given this notation, define the following $n$-state nonlinear system:
\begin{align} \label{eq:sys}
    \dot{x}(t) = c + Lx(t) + f(x(t)),
\end{align}
where $x(t) \in \Rn{n}$ is the state vector at time $t$. The time dependency of the state vector $x(t)$ is occasionally omitted to simplify the notation. Furthermore, we assume the quadratic nonlinearity is lossless: $x^\top f(x) = 0$
for all $x\in \Rn{n}$. The lossless condition is: 
\begin{align}
    0 = x^\top f(x) =  \sum_{i=1}^n x_i \left( x^\top Q^{(i)} x \right) = \sum_{i,j,k=1}^n Q_{jk}^{(i)} x_i x_j x_k \qquad  \forall x \in \Rn{n}.
\end{align}
Note that every term in the product $x^\top f(x)$ is cubic in the entries of $x$. The lossless condition holds if and only if the coefficient of each term $x_i x_j x_k$ is zero. Accounting for the fact that $x_i x_j x_k = x_j x_k x_i = x_k x_i x_j$, the nonlinearity $f(x)$ is energy preserving if and only if the symmetric matrices $\{Q^{(i)}\}_{i=1}^n$ satisfy the following condition:
\begin{align}   \label{eq:lossless}
        Q^{(i)}_{jk} + Q^{(j)}_{ik} + Q^{(k)}_{ij} = 0, & &\forall i,j,k 
        \in \{ 1,\dots,n \}.
\end{align}

One important feature of a quadratic system is that the nonlinearities are preserved under coordinate translation. Specifically, define $y:=x-m$, where $m\in\Rn{n}$ is a constant translation. The dynamics after translation are:
\begin{align} \label{eq:sys_shifted}
    \dot{y}(t) &= d(m) + A(m)\, y(t) + f(y(t)), 
\end{align}
where $d(m) := c + Lm + f(m)$ and 
\begin{align}
    A(m) := L + 2\bmtx m^\top Q^{(1)} \\ \vdots \\ m^\top Q^{(n)} \emtx.
\end{align}
The quadratic nonlinearity in~\eqref{eq:sys_shifted} is exactly the same as in~\eqref{eq:sys}. Hence the nonlinearity in~\eqref{eq:sys_shifted} is also lossless. The symmetric part of $A(m)$ is given by $A_s(m):=\frac{1}{2}(A(m) + A(m)^\top )$.
This symmetric part plays an important role in the boundedness analysis introduced later. It can be expressed simply as:
\begin{align}
    \label{eq:As}
    A_s(m) 
    = L_s - \sum_{k=1}^n m_k Q^{(k)},
\end{align}
where $L_s=\frac{1}{2}(L+L^\top)$ is the symmetric part of $L$ and $m_k$ is the $k^{th}$ element of the coordinate shift $m$. Equation~\ref{eq:As} is due to the lossless property of the nonlinearity~\eqref{eq:lossless}:
\begin{align}
    \bmtx m^\top Q^{(1)} \\ \vdots \\ m^\top Q^{(n)} \emtx_{i,j} + \bmtx Q^{(1)}m & \ldots & Q^{(n)}m \emtx_{i,j} = \sum_{k=1}^n m_k Q_{jk}^{(i)} + \sum_{k=1}^n m_k Q_{ik}^{(j)} = - \sum_{k=1}^n m_k Q_{ij}^{(k)}.
\end{align}

\subsection{Boundedness of a System}

Let $\phi(t; x_0, t_0)$ denote the solution of~\eqref{eq:sys} at time $t\geq t_0$ starting from the initial condition $x_0 \in \Rn{n}$ at time $t_0$. As the stability of a nonlinear system is hard to analyze in general, the boundedness of solutions $\phi(t; x_0, t_0)$ serves as an alternative notion to characterize the behavior of a system. The boundedness of the system can be defined as in Chapter 4.8 of Khalil~\cite{khalil2002nonlinear}.
\begin{definition}[Boundedness]     \label{def:Boundedness}
    The solutions of system~\eqref{eq:sys} are globally uniformly ultimately bounded if there exists a scalar $\beta>0$, independent of $t_0\geq 0$, and a function $T:\Rn{} \to \Rn{}$ such that: 
    \begin{align}
        \norm{\phi(t; x_0, t_0)} \leq \beta & & \forall  x_0\in\Rn{n} \text{ and } \forall t \geq t_0 + T(\norm{x_0}).
    \end{align}
\end{definition}
In words, a system is globally uniformly ultimately bounded with bound $\beta > 0$ if the trajectory from any initial condition eventually converges to within a spherical ball centered at the origin and with radius $\beta$. Moreover, $T(\norm{x_0})$ is the time required to converge within the bound $\beta$. This time is independent of $t_0$ but depends on the norm of the initial condition. In this paper, only globally uniformly ultimate boundedness is investigated and we will simply refer to this as the ``boundedness" in the remainder of the paper. 

Note that the boundedness of the shifted system~\eqref{eq:sys_shifted} is equivalent to the boundedness of the original system~\eqref{eq:sys} as the coordinates only differ by translation. 
In addition, the triangle inequality implies 
$\norm{x(t)} \le \norm{y(t)} + \norm{m}$. Hence,  if a trajectory in shifted coordinates is bounded by $\beta$, i.e. $\norm{y(t)}\le\beta$, then the corresponding trajectory in the original coordinates is bounded by $\beta + \norm{m}$. Hence, one can show the boundedness of the original system~\eqref{eq:sys} by showing the boundedness of the shifted system~\eqref{eq:sys_shifted}. 


\subsection{A Sufficient Condition for Boundedness} \label{sect:prel_LorenzCond}

This section reviews a sufficient condition to certify the boundedness of quadratic systems~\eqref{eq:sys}. This section is mainly based on early results discussed by Lorenz~\cite{lorenz1963deterministic}, which were further formalized by Schlegel and Noack~\cite{schlegel2015long}. First, we introduce the notion of a monotonically attracting trapping region as defined in Schlegel and Noack~\cite{schlegel2015long}. 
\begin{definition}[Trapping region]
A trapping region $D \subseteq \Rn{n}$ is a compact set that is forward invariant, i.e., if $x(t_0) \in D$ then $x(t) \in D$ for all $t \ge t_0$. In other words, once a trajectory enters a trapping region, it remains in the trapping region for all future time. A trapping region is termed globally monotonically attracting if an energy function is strictly monotonically decreasing along all trajectories starting from an arbitrary state outside of $D$.
\end{definition}
The existence of a monotonically attracting trapping region for a system implies that the system is bounded, as all trajectories will enter the compact set $D$. We will focus on globally monotonically attracting trapping regions and will refer to this simply as a ``trapping region." 

Section 2 of Lorenz's 1963 paper~\cite{lorenz1963deterministic} provides a condition for the existence of a trapping region for the dynamical system~\eqref{eq:sys}. The condition uses the kinetic energy $K_0(x) := \frac{1}{2}\norm{x}^2$. The derivative of this energy function along trajectories of the system~\eqref{eq:sys} is:
\begin{align}
    \label{eq:dK0dt}
    \frac{d}{dt}K_0\left(x(t)\right) := \frac{1}{2} \left( x(t)^\top\dot{x}(t)
      + \dot{x}(t)^\top x(t) \right)= c^\top x(t) + x(t)^\top L_s x(t).
\end{align}
The derivative of $K_0(x(t))$ does not depend on the quadratic dynamics $f$ as the nonlinearity is lossless. Moreover, the quadratic term in \eqref{eq:dK0dt} dominates the linear term if $\norm{x}$ is sufficiently large and $L_s$ is non-singular. Thus, the energy of all trajectories decreases if $L_s$ is negative definite and $\norm{x}$ is large enough. The energy can increase only when close to the origin. Hence, a trapping region around the origin exists as all trajectories converge toward the origin or stay close to the origin. In summary, if $L_s$ is negative definite, a trapping region exists, and the system is bounded.

In fact, a more precise characterization for a trapping region can be obtained using the eigenvalues of $L_s$. Assume $L_s$ is negative definite and denote the eigenvalues by $\lambda_n \le \dots \le \lambda_1 <0$. Define the closed-ball with center $v$ and radius $R$ as $B(v,R):=\{ x\in\Rn{n} \, : \, \norm{x-v} \le R\}$. The time-derivative of $K_0$ is bounded using the largest (most positive) eigenvalue and the Cauchy-Schwartz inequality:
\begin{align} \label{eq:estimating_R0}
    \frac{d}{dt}K_0(x(t)) \le
        \norm{ c } \, \norm{x(t)}
        + \lambda_1 \norm{x(t)}^2
        = \left( \norm{c} + \lambda_1 \norm{x(t)} \right) \cdot \norm{x(t)}.
\end{align}
Thus the derivative of $K_0$ is negative if $\norm{x(t)} > \frac{\norm{c}}{|\lambda_1|}$. Hence $B(0,R_0)$ is a trapping region for the system~\eqref{eq:sys} in original coordinates with the radius
$R_0:=\frac{\norm{c}}{|\lambda_1|}$. Note that this bound is conservative in general and is attained only when the vector $c$ aligns exactly with the eigenvector corresponding to $\lambda_1$.

Next, define the set $E:=\left\{ x(t)\in\Rn{n}  \, : \, \frac{d}{dt} K_0(x(t)) \ge 0 \right\}$. The set $E$ corresponds to states where the energy is (instantaneously) non-decreasing. $B(0,R_0)$ is shown above to be a trapping region, and the energy is strictly decreasing outside this ball. Hence $E \subseteq B(0,R_0)$. To further characterize $E$, assume $L_s$ is negative definite and use completion of squares to rewrite  \eqref{eq:dK0dt} as follows:
\begin{align} \label{eq:dK0dt_Ellipse}
    \frac{d}{dt} K_0(x(t)) = 
     \left( x(t) + \frac{1}{2} L_s^{-1} c \right)^\top L_s
     \left( x(t) + \frac{1}{2} L_s^{-1} c \right)
      - \frac{1}{4} c^\top L_s^{-1} c.
\end{align}
The boundary of $E$ corresponds to the states where
$\frac{d}{dt} K_0(x(t)) = 0$. It follows from
\eqref{eq:dK0dt_Ellipse} that this boundary is
an ellipsoid.\footnote{The matrix $L_s$ is symmetric
and negative definite. Thus it has an eigenvalue
decomposition $L_s=U \Lambda U^\top$ where $U$ is orthogonal and $\Lambda<0$ is a diagonal matrix with the eigenvalues  $\lambda_n \le \dots \le \lambda_1 <0$ along the diagonal. Define the rotated coordinates $z(t):=U^\top x(t)$. The condition $\frac{d}{dt} K_0(x(t)) = 0$ can be written in these rotated coordinates as
$(z(t)-v)^\top \, \Lambda \, (z(t)-v) =\frac{1}{4}c^\top L_s^{-1}c$ where $v:=-\frac{1}{2}  U^\top L_s^{-1}  c \in \Rn{n}$
is the ellipsoid center in the rotated coordinates.  This can be further simplified to 
$\sum_{i=1}^n \frac{1}{\alpha_i^2} ( z_i- v_i)^2 = 1$ where $
 \alpha_i := \frac{1}{2} \sqrt{ \frac{ c^\top L_s^{-1} c } {\lambda_i} }$ are the semi-axis lengths.} 
The point $x(t)=-\frac{1}{2}L_s^{-1}c \in E$ is the ellipsoid center. This is the point of maximum energy growth: $\frac{d}{dt} K_0(x(t)) = - \frac{1}{4} c^\top L_s^{-1} c >0$ (assuming $c\ne 0$).

Note that the system~\eqref{eq:sys} can have more than one equilibrium point. The analysis given above is independent of the equilibrium point and only depends on $c$ and $L_s$. This allows analyzing the boundedness of systems with multiple equilibrium points or even chaotic behavior. In fact, if a trapping region exists, all equilibrium points would be inside the trapping region and on the boundary of the non-decreasing energy set $E$. This is because the energy evolution at an equilibrium is always $0$. 

\begin{figure}
    \centering
    \includegraphics[width=0.6\linewidth]{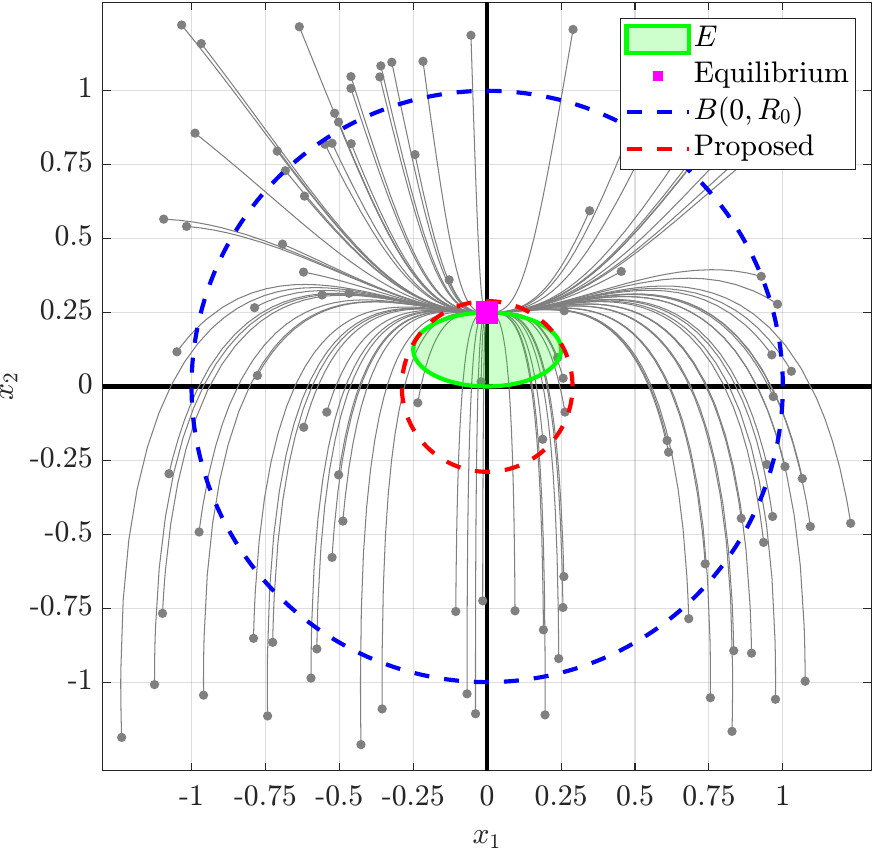}
    \caption{An illustration of trapping region for the two-state system~\eqref{eq:TwoStateSys} with kinetic energy $K_0(x)=\frac{1}{2}x^\top x$. All trajectories with random initial conditions (gray) converge to the equilibrium point at $(0,0.25)$ (magenta square). The set of states with non-decreasing energy $E$ (green solid ellipsoid) and a conservative trapping region $B(0,R_0)$ (blue-dashed circle) are shown. The tightest trapping region centered at the origin can be computed by our proposed method (Section~\ref{sect:Verif_Computing_beta}) and is shown in the red-dashed circle.}
    \label{fig:TwoStateSys_Prelim}
\end{figure}

These theoretical results are illustrated in Figure~\ref{fig:TwoStateSys_Prelim} for the following two-state system with an lossless quadratic term:
\begin{align}   \label{eq:TwoStateSys}
    \frac{d}{dt} \bmtx x_1 \\ x_2 \emtx = 
    \bmtx 0 \\ 1 \emtx + \bmtx -1 & 0 \\ 0 & -4 \emtx x +
    \bmtx -x_1x_2 \\ x_1^2 \emtx. 
\end{align}
The trajectories (light gray lines) in Figure~\ref{fig:TwoStateSys_Prelim} are generated with random initial conditions and converge to the single stable equilibrium $(0,0.25)$. The constant vector in the dynamics is $c=\bsmtx 0 \\ 1 \esmtx$ and the largest (most positive) eigenvalue of $L_s$ is $\lambda_1 = -1$. Hence, $B(0,R_0)$ is a trapping region with $R_0=\frac{\norm{c}}{|\lambda_1|}=1$ as derived above. This is shown as the blue-dashed circle in Figure~\ref{fig:TwoStateSys_Prelim}. All trajectories eventually converge inside this ball as expected. The green ellipsoid in this figure corresponds to the states where $\frac{d}{dt}K_0(x(t)) \geq 0$. Based on~\eqref{eq:dK0dt_Ellipse}, the boundary of this ellipse is defined by:
\begin{align*}
    \frac{d}{dt}K_0(x(t)) 
    &= -x_1^2 - 4\left(x_2 - \frac{1}{8}\right)^2 + \frac{1}{16} = 0.
\end{align*}
As shown in Figure~\ref{fig:TwoStateSys_Prelim}, the circle $B(0,R_0)$ is a conservative estimate as any ball centered at the origin containing the ellipsoid $E$ is a trapping region. This conservatism is rooted in the derivation using Cauchy-Schwartz inequality in~\eqref{eq:estimating_R0}. We will derive in Section~\ref{sect:Verif_Computing_beta} a convex optimization to compute the tightest ball containing the ellipsoid $E$. This tightest (least conservative) ball is shown as the red-dashed circle in Figure~\ref{fig:TwoStateSys_Prelim}.

\subsection{Boundedness Using Coordinate Translations}
\label{sect:Prel_SchlegelNoack}

The previous subsection presented a sufficient condition for boundedness using the kinetic energy centered at the state origin. Schlegel and Noack~\cite{schlegel2015long} further generalized the condition using coordinate translations. Consider the translated system~\eqref{eq:sys_shifted} with $y(t):=x(t)-m$. The kinetic energy for the translated system
is $K_m(y) := \frac{1}{2}\norm{y}^2$ with corresponding
time derivative:
\begin{align}
   \label{eq:dKmdt}
    \frac{d}{dt}K_m(y(t)) 
    := \frac{1}{2} \left( y(t)^\top\dot{y}(t) + \dot{y}(t)^\top y(t)\right) = d(m)^\top y(t) + y(t)^\top A_s(m) y(t).
\end{align}
Using the same argument as in Section~\ref{sect:prel_LorenzCond}, if the matrix $A_s(m)$ is negative definite, then
the shifted system~\eqref{eq:sys_shifted} has a trapping region $B(0,R_m)$ in the shifted coordinate for some $R_m$.  This implies the original system~\eqref{eq:sys} is bounded by $\norm{m}+R_m$. Furthermore, a valid radius is
$R_m:=\frac{\norm{d(m)}}{|\lambda_1|}$ where $\lambda_1$
is the largest (most positive) eigenvalue of $A_s(m)$.
The following theorem in Schlegel and Noack~\cite{schlegel2015long} formalizes these discussions for monotonically attracting trapping regions.

\begin{theorem}[Theorem 1 of Schlegel and Noack{~\cite{schlegel2015long}}] \label{theorem:TrappingRegion}
    The system~\eqref{eq:sys} has a monotonically attracting trapping region $B(m, R_m)$ if and only if there exists an $m\in\Rn{n}$ such that the real, symmetric matrix $A_s(m)$ is negative definite. If $A_s(m)$ has eigenvalues $\lambda_n \leq \dots \leq \lambda_{1}<0$, then a trapping region is given by $B(m,R_m)$ with the radius $R_m:=\frac{\norm{d(m)}}{|\lambda_1 |}$.
\end{theorem}
A formal proof of this necessary and sufficient trapping region condition is given in Appendix B of Schlegel and Noack~\cite{schlegel2015long}. Note that the radius $R_m$ given in Theorem~\ref{theorem:TrappingRegion} is different from the statement in Schlegel and Noack~\cite{schlegel2015long} as there is a minor typo that is also pointed out in Kaptanoglu et al~\cite{Kaptanoglu2021}. Also, the proof in Schlegel and Noack~\cite{schlegel2015long} regarding the estimated radius $R_m$ used a geometrical argument instead of the Cauchy-Schwartz inequality approach presented in~\eqref{eq:estimating_R0}. The two estimates are the same, and the Cauchy-Schwartz inequality is adopted in this paper for simplicity.

The key point of Theorem~\ref{theorem:TrappingRegion} is that the boundedness of the system~\eqref{eq:sys} can be determined by an algebraic condition. In other words, if $m$ exists such that $A_s(m)$ is negative definite, then the system has a trapping region $B(m,R_m)$. This corresponds to a bound $R_m$ in the shifted coordinates and $\norm{m}+R_m$ in the original coordinates.  This approach does not require solving the nonlinear differential equation for all possible initial conditions. Note that the converse is not true, i.e., a bounded system does not imply the existence of a trapping region. Consider a trivial system $\dot{x}=0$. All solutions are stationary and thus bounded. However, this system does not have a monotonically attracting trapping region, i.e., one cannot find an $m$ such that $A_s(m)$ is negative definite. If no trapping region exists, then Schlegel and Noack~\cite{schlegel2015long} point out that additional analysis of the linear and nonlinear terms is needed to determine whether the system is bounded. 

Even though Theorem~\ref{theorem:TrappingRegion} provides an approach to verify the boundedness of the system, there are two remaining caveats. First, verifying this algebraic condition is non-trivial as it involves a search over $m\in\Rn{n}$. An inf-sup problem is formulated in Schlegel and Noack~\cite{schlegel2015long} to minimize the largest eigenvalue of $A_s(m)$ and simulated annealing is proposed as a solution method. However, this formulation is a non-convex optimization, and hence there are no guarantees of convergence to the global optimum. Thus the test is inconclusive if the simulated annealing fails to find a feasible $m$ with the largest eigenvalue of $A_s(m)$ being negative. In other words, a trapping region may exist even if the simulated annealing algorithm fails. The second caveat of Theorem~\ref{theorem:TrappingRegion} is that the estimated bound $R_m$ is an approximate upper bound and is usually conservative. As mentioned in Section~\ref{sect:prel_LorenzCond}, the approximation is made by the analysis of eigenvalues of symmetric matrix $A_s(m)$ with the worst-case alignment of direction $d(m)$ in the Cauchy-Schwartz inequality~\eqref{eq:estimating_R0}. Hence, the upper bound approximation $R_m$ is often conservative. We address both of these caveats in the next section by providing: (i) a convex optimization to compute the trapping region, and (ii) a tight bound on the trapping region radius.

\section{Computation of Trapping Regions Via Convex Optimization}
\label{sect:Verif}

This section revisits the boundedness condition in Theorem~\ref{theorem:TrappingRegion}. First, the condition is reformulated as a constraint in the form of a linear matrix inequality (LMI), whose feasible set is a convex set. The feasibility of this LMI condition can be tested numerically, and a feasible solution, if one exists, verifies the existence of a trapping region. This avoids the non-convex inf-sup formulation in Schlegel and Noack~\cite{schlegel2015long}. Second, a quadratically constrained quadratic program (QCQP) is derived to compute the exact radius of the bounded region without conservative approximation. Lastly, the set of furthest points with non-decreasing energy can be identified by the dual solution of the QCQP and shown to be a hypersphere. 

\subsection{Trapping Condition as a Convex Constraint} \label{sect:Verif_Boundedness}

By Theorem~\ref{theorem:TrappingRegion}, the system~\eqref{eq:sys} has a trapping region if and only if there exists $m\in\Rn{n}$ such that $A_s(m)$ is negative definite, i.e.,
\begin{align}   \label{eq:LMI_boundedness}
    A_s(m) = L_s - \sum_{i=1}^N m_i Q^{(i)} \prec 0.  
\end{align}
%
Observe that $A_s(m)$ is an affine function of $m$. Hence, the negative definite constraint~\eqref{eq:LMI_boundedness} is a linear matrix inequality (LMI) in the variable $m$~\cite{boyd2004convex, boyd1994linear}. LMIs have been widely studied in the optimization community and are known to have a convex feasible set, i.e., a convex combination of two solutions is also a solution (Chapter 2.3.2 of Boyd and Vandenberghe~\cite{boyd2004convex}). Hence, the set of $m\in\Rn{n}$ that satisfy~\eqref{eq:LMI_boundedness} is a convex set. LMI feasibility problems can be efficiently solved with many numerical optimization solvers, e.g., SeDuMi~\cite{SeDuMi}, MOSEK~\cite{mosek}, and SDPT3~\cite{tutuncu2003SDPT3}. There is also software that enables easy implementation and interface with these solvers, e.g., CVX~\cite{grant2014cvx}. Finally, LMI feasibility problems have useful duality results such as the theorem of alternatives (Chapter 5.8 of Boyd and Vandenberghe~\cite{boyd2004convex}). Thus, the numerical algorithm will return a coordinate shift $m\in\Rn{n}$ satisfying~\eqref{eq:LMI_boundedness} (and hence verifying the existence of a trapping region) or it will return a dual variable that certifies that no such feasible $m$ exists.

The LMI feasibility problem can be re-formulated as an optimization to minimize the largest eigenvalue of $A_s(m)$:
\begin{mini}
    {m\in\Rn{n}, a\in\Rn{}}{a}{\label{SDP:min_lambda}}{a^*=}
    \addConstraint{A_s(m) \preceq a I_n}.
\end{mini}
This optimization has a single LMI constraint with the objective given by the slack variable $a$. This falls into the category of Semidefinite Programs (SDPs) (Chapter 4.6.2 of Boyd and Vandenberghe~\cite{boyd2004convex}). SDPs are convex optimizations that can also be efficiently solved using the numerical tools mentioned above. The particular problem in~\eqref{SDP:min_lambda} can be solved on a standard desktop for problems with dimensions up to $n \sim O(100)$. Larger problems require more specialized implementations. 
In contrast to the inf-sup formulation proposed in Schlegel and Noack~\cite{schlegel2015long}, this SDP formulation is guaranteed to converge to the global optimum (within numerical tolerances). Also, the SDP~\eqref{SDP:min_lambda} is always feasible by choosing any $m\in\Rn{n}$ and a sufficiently large value of $a$. If the optimal cost satisfies $a^* < 0$, then there exists $m$ such that $A_s(m)$ is negative definite. In this case, a trapping region exists, and the system~\eqref{eq:sys} is bounded. If $a^* \geq 0$, then there is no feasible $m$, and no trapping region exists. 


\subsection{Computing the Tightest Boundedness Region with SDP} \label{sect:Verif_Computing_beta}

Assume the SDP~\eqref{SDP:min_lambda} has an optimal cost $a^*<0$ so that
there exists an $m$ such that $A_s(m)\prec 0$.  By Theorem~\ref{theorem:TrappingRegion}, the system~\eqref{eq:sys} has a monotonically attracting trapping region. Let $\lambda_1$ denote the largest (least negative) eigenvalue of $A_s(m)$. A conservative bound on the trapping region is given by $B(m,R_{m})$ with the radius $R_m:=\frac{\norm{d(m)}}{|\lambda_1 |}$. In this subsection, we utilize duality in optimization theory to compute the smallest (least conservative) radius $R_m^*$ such that $B(m,R_m^*)$ is a trapping region.

If $\frac{d}{dt} K_m(y(t)) \ge 0$ then the energy is instantaneously non-decreasing. It follows from \eqref{eq:dKmdt} that
the set of all states where the energy is instantaneously non-decreasing (in shifted coordinates) is: 
\begin{align}   \label{eq:EnergyNonDecreasingSet}
     E:=\left\{ y \in\Rn{n}  \, : \, d(m)^\top y + y^\top A_s(m) y \ge 0 \right\}.
\end{align}  
The smallest (least conservative) radius $R_m^*$ corresponds to the smallest ball (in shifted coordinates) that contains all states in $E$. Thus
the smallest radius is given by:
\begin{align*}
    R_m^* :=  \min_{R_m} & \quad R_m   \\ 
            \text{s.t.}  & \quad E \subseteq B(0,R_m).
\end{align*}
Equivalently, we can find the state in $E$ with the largest norm, i.e., maximize $\norm{y}^2$ subject to $y\in E$. This equivalent formulation leads to the following optimization:
\begin{subequations}  \label{eq:QCQP}
    \begin{align}
      (R_m^*)^2 = \max_y \quad &f_0(y):= y^\top y \\
      \text{s.t. } \quad &f_1(y):=  d(m)^\top y + y^\top A_s(m) y \geq 0.     
    \end{align}
\end{subequations}
Note that \eqref{eq:QCQP} is a non-convex problem as we are maximizing a convex objective. Non-convex problems can, in general, have multiple local optima and solving for global optima can be difficult.  In fact, \eqref{eq:QCQP} is a quadratically constrained quadratic program (QCQP) with only one quadratic constraint. This particular problem, while nonconvex, can be efficiently solved for the global optimum. The Lagrange dual of the QCQP~\eqref{eq:QCQP} (Appendix B of Boyd and Vandenberghe~\cite{boyd2004convex}) is:
\begin{subequations} \label{eq:QCQPdual}
\begin{align} 
    (R_m^*)^2 = \min_{\lambda\geq 0, \gamma} \quad &\gamma \quad \\
    \text{s.t. }  \quad &f_2(\lambda, \gamma) :=    \bmtx I+\lambda A_s(m) & \frac{\lambda}{2}d(m) \\ \frac{\lambda}{2}d(m)^\top & -\gamma \emtx \preceq 0.
\end{align}
\end{subequations}
The optimal cost for this dual SDP~\eqref{eq:QCQPdual} is identical to the optimal cost of the original QCQP~\eqref{eq:QCQP}. This follows from strong duality if $d(m)\ne 0$.  In particular, the point $y_0=-\frac{1}{2}\frac{d(m)^\top d(m)}{d(m)^\top A_s(m) d(m)}\, d(m)$ is strictly feasible for the QCQP~\eqref{eq:QCQP}, i.e. $f_1(y_0)<0$.  Hence the QCQP satisfies Slater’s constraint qualification and strong duality holds (Appendix B of Boyd and Vandenberghe~\cite{boyd2004convex}).

If $d(m)=0$, then the constraint $f_1(y) \leq 0$ does not have a constant term and is not strictly feasible.  The only feasible point of the QCQP~\eqref{eq:QCQP} in this case is $y_0=0$ and hence the optimal cost is $R_m^*=0$. The dual~\eqref{eq:QCQPdual} has the same optimal cost of 0 when $d(m)=0$. This optimal dual cost is obtained by $\gamma=0$ with $\lambda>0$ sufficiently large such that $I+\lambda A_s(m) \preceq 0$. This is possible 
since $A_s(m)\prec 0$. Thus the QCQP and dual achieve the same cost even when $d(m)=0$.  The optimal cost $R_m^*=0$ implies that the smallest boundedness region is the origin in the shifted coordinates.  Thus $A_s(m)\prec 0$ and $d(m)=0$ implies that the state $m$ is a globally stable equilibrium point for the system~\eqref{eq:sys}.

In summary, we can characterize, with no conservatism, the smallest radius $R_m^*$ such that the trapping region $B(m, R_m^*)$ contains all the states where the energy is instantaneously non-decreasing.  All trajectories of the original system~\eqref{eq:sys} will eventually converge to and be tightly bounded in the region $B(m, R_m^*)$. If $d(m)=0$, then $R_m^*=0$ and the state $m$ corresponds to a globally stable equilibrium. For the case $d(m)\neq 0$, the dual SDP~\eqref{eq:QCQPdual} can be solved numerically to obtain optimal radius $R_m^*$. 
The global optimum $R_m^*$ from the the QCQP~\eqref{eq:QCQP} and its dual SDP~\eqref{eq:QCQPdual} provides a tighter bound than  the value $R_m:=\frac{\norm{d(m)}}{|\lambda_1 |}$ given Theorem~\ref{theorem:TrappingRegion}. 
The discussion in this section applies to any $m$ such that $A_s(m)\prec 0$.  It can also be applied to the optimal value $m^*$ obtained from the SDP~\eqref{SDP:min_lambda} assuming it has an optimal cost $a^*<0$.


Note that one might potentially want to find the tightest trapping region among all possible coordinate shifts. This corresponds to jointly optimizing the problem~\eqref{eq:QCQP} over $(m,y)$ with the additional constraint~\eqref{eq:LMI_boundedness} such that a trapping region exists. Recall that $d(m) = c + Lm + f(m)$.  Hence constraints of~\eqref{eq:QCQP} have a term $f(m)^\top y$ that is cubic in the variables $(m,y)$. The problem no longer has the same strong duality results, in general. We are not aware of an efficient approach to jointly optimize the problem~\eqref{eq:QCQP} over $(m, y)$. Nonlinear optimization algorithms can be applied but they are not guaranteed to converge to the global optima. Further investigation of solving this problem is left to future work.

\subsection{Recovering Furthest Points with Non-Decreasing Energy} \label{sect:Verif_ystar}

In this subsection, we characterize the set of optimal solutions $y^*$ of the QCQP~\eqref{eq:QCQP}. These correspond to states $y^*$ with norm of $R_m^*$ and non-decreasing energy. The case $d(m)=0$ has the trivial optimal solution $y^*=0$ as discussed in the previous subsection. This subsection focuses on the case of $d(m)\neq 0$.

Strong duality holds between the QCQP~\eqref{eq:QCQP} and the dual SDP~\eqref{eq:QCQPdual} when $d(m) \neq 0$. Hence any primal and dual optimal solutions must satisfy the first-order optimality conditions. Specifically, any primal and dual optimal variables $(y^*,\gamma^*,\lambda^*)$ must satisfy the following  Karush–Kuhn–Tucker (KKT) conditions (Chapter 5.5.3 of Boyd and Vandenberghe~\cite{boyd2004convex}):
\begin{subequations}    \label{QCQP_KKT}
    \begin{align}
        &\textit{Stationary condition:} &&   \left.\nabla_y \left[f_0(y) + \lambda^* f_1(y) \right] \right|_{y=y^*} = 0      \label{QCQP_KKT:Stationary}    \\  
        &\textit{Primal feasibility:} &&    f_1(y^*) \leq 0          \label{QCQP_KKT:Primal}       \\  
        &\textit{Dual feasibility:}   &&  \lambda^* \geq 0, \quad f_2(\lambda^*,\gamma^*) \succeq 0                  \label{QCQP_KKT:Dual}         \\  
        &\textit{Complementary slackness:} && \lambda^* f_1(y^*) = 0                \label{QCQP_KKT:CompSlack}               
    \end{align}
\end{subequations}

The stationary condition~\eqref{QCQP_KKT:Stationary} for the given $f_0$ and $f_1$ is:
\begin{align*}   
    0 = \left.\nabla_y \left[f_0(y) + \lambda^* f_1(y) \right] \right|_{y=y^*} &= 2(I + \lambda^* A_s(m)) y^* + \lambda^* d(m).
\end{align*}
Let $r$ denote the rank of $(I+\lambda^* A_s(m))$ and $V\in \Rn{n\times (n-r)}$ be a matrix whose columns are orthonormal and span the null space of $(I+\lambda^* A_s(m))$.
Then for any optimal point $y^*$ there exists $c \in \Rn{n-r}$ such that:
\begin{align*}
    y^*(c) &= -\frac{\lambda^*}{2}(I+\lambda^* A_s(m))^\dagger d(m) + Vc,
\end{align*}
where $(\cdot)^\dagger$ is the pseudoinverse of a matrix.  

Next, observe that the dual constraint $f_2(\lambda^*, \gamma^*) \succeq 0$ in~\eqref{QCQP_KKT:Dual} is infeasible when $\lambda^* = 0$. Hence, the dual constraints~\eqref{QCQP_KKT:Dual} imply that $\lambda^*>0$. The complementary slackness condition~\eqref{QCQP_KKT:CompSlack} then further implies that $f_1(y^*(c)) = 0$. Note that $f_1(y^*(c)) = 0$ also implies the primal feasibility~\eqref{QCQP_KKT:Primal}. Therefore, the primal-dual point $(y^*(c), \lambda^*, \gamma^*)$ satisfies KKT conditions~\eqref{QCQP_KKT} for any $c\in\Rn{n-r}$. The primal points $y^*(c)$ are only candidates for the optimal solution of the QCQP~\eqref{eq:QCQP} as the KKT conditions are only necessary for optimality.  The null space vector $c^*$ yields an optimal point if $y^*(c^*)$ achieves the optimal cost:
\begin{align}
    (R_m^*)^2 &= f_0(y^*(c^*)) = \frac{(\lambda^*)^2}{4}d(m)^\top \left( (I+\lambda^* A_s(m))^\dagger \right)^2 d(m) + (c^*)^\top c^* .
\end{align}
Therefore, the set of optimal solutions $y^*(c^*)$ forms a $(n-r)$-dimensional hypersphere with radius of $\norm{c^*}$:
\begin{align}   \label{eq:critical_ystar}
    \left\{ y^*(c^*) =  -\frac{\lambda^*}{2}(I+\lambda^* A_s(m))^\dagger d(m) + Vc^*
      \, : \, \norm{c^*}^2 = (R_m^*)^2 - \frac{(\lambda^*)^2}{4}d(m)^\top \left( (I+\lambda^* A_s(m))^\dagger \right)^2 d(m) \right\},
\end{align}
where the first term of $y^*(c^*)$ is the center of the hypersphere, columns of $V$ are $(n-r)$ directions of the hypersphere in $n$-dimensional space. This $(n-r)$-dimensional hypersphere corresponds to the points $y^*$ furthest from the origin of the shifted system~\eqref{eq:sys_shifted} where the energy $K_m(y^*(c^*))$ is instantaneously non-decreasing.

\subsection{Comparison to Existing Results}
\label{sect:comparison}

The following theorem summarizes the proposed convex optimization method for computing a trapping region.
\begin{theorem}~\label{theorem:TrappingRegion_SDP}
    The system~\eqref{eq:sys} has a monotonically attracting trapping region $B(m, R_m^*)$ if and only if 
    the optimal value of SDP~\eqref{SDP:min_lambda} satisfies $a^*<0$. 
    Moreover, if $a^*<0$ then there exists at least one $m$ such that $A_s(m) \prec 0$.  The tightest trapping region centered at any such $m$ is given by $B(m, R_m^*)$ where $(R_m^*)^2$ is the optimal cost of~\eqref{eq:QCQP}. 
    
\end{theorem}
The first statement of Theorem~\ref{theorem:TrappingRegion_SDP} follows from Theorem~\ref{theorem:TrappingRegion} and the definition of the SDP~\eqref{SDP:min_lambda}. Specifically, there exists $m\in\Rn{n}$ such that $A_s(m)$ is negative definite if and only if the SDP~\eqref{SDP:min_lambda} has optimal cost satisfying $a^*<0$. Theorem~\ref{theorem:TrappingRegion} only provides an upper bound on the trapping region radius. The second statement of Theorem~\eqref{theorem:TrappingRegion_SDP} is that the tightest trapping region is given by the optimization~\eqref{eq:QCQP}. As discussed in Section~\ref{sect:Verif_Computing_beta}, the optimization~\eqref{eq:QCQP} can be equivalently formulated as a dual SDP~\eqref{eq:QCQPdual}. In summary, the numerical method is to first solve the SDP~\eqref{SDP:min_lambda} and, if $a*<0$, then solve for $R_m^*$ using the SDP in~\eqref{eq:QCQPdual}. As an aside, the SDP~\eqref{eq:QCQPdual} also yields the optimal Lagrange multiplier $\lambda^*$.  This can be used to compute the points on the boundary of $B(m,R_m^*)$ with non-decreasing energy. These correspond to the points $x^* = y^*(c) +m$ with $y*$ in the set~\eqref{eq:critical_ystar}.



Finally, this paper has three key differences compared to the result in Schlegel and Noack~\cite{schlegel2015long}. First, this paper provides a computational approach to verify the existence of a trapping region with guarantees, while Schlegel and Noack's non-convex approach does not provide such guarantees. Second, the tightest trapping region is identified in this paper instead of an upper bounding radius in Schlegel and Noack. Third, this paper characterizes the set of points on the boundary of the tightest trapping region with non-decreasing energy, which is not discussed in Schlegel and Noack.

\section{Numerical Examples}    \label{sect:Examp}


Here, we provide numerical examples to illustrate the proposed SDP-based analysis. First, the trapping region of the two-dimensional example system~\eqref{eq:TwoStateSys} is revisited. Second, the Lorenz attractor is presented as a benchmark to compare our method against that given in Schlegel and Noack~\cite{schlegel2015long}. Lastly, a system with an adjustable number of states is used to study the scalability of the proposed method. All analyses are implemented in Matlab with CVX~\cite{grant2014cvx} and SDP solver MOSEK~\cite{mosek}. The examples are run on a standard desktop with an Intel(R) i7-9700 CPU @ 3.00GHz and 16 GB memory. Complete source code for the implementation and examples can be found on GitHub at \url{https://github.com/SCLiao47/Boundedness_LosslessQuadSys}.

\subsection{Bounded Two-Dimensional System}
Recall the two-dimensional system previously introduced in Section~\ref{sect:prel_LorenzCond}:
\begin{align}   \label{eq:TwoStateSys_example}
    \frac{d}{dt} \bmtx x_1 \\ x_2 \emtx = 
    \bmtx 0 \\ 1 \emtx + \bmtx -1 & 0 \\ 0 & -4 \emtx x +
    \bmtx -x_1x_2 \\ x_1^2 \emtx. 
\end{align}
The linear symmetric part is $L_s = \bsmtx -1 & 0 \\ 0 & -4 \esmtx$. If there is zero coordinate shift, i.e., $m=0$, then $A_s(m)=L_s$ is negative definite. Hence the system is bounded by Theorem~\ref{theorem:TrappingRegion}, and a conservative estimate on the trapping region radius was given as $R_0=\frac{\norm{c}}{|\lambda_1|} = 1$. The optimal cost for the SDP~\eqref{SDP:min_lambda} satisfies $a^*<0$ because $A_s(m) \prec 0$ at $m=0$. Moreover, we can use Theorem~\ref{theorem:TrappingRegion_SDP} to compute the tightest trapping region $B(0, R_0^*)$. The optimal solution of the SDP~\eqref{eq:QCQPdual} is $R_0^*=0.289$ with $\lambda^*=1$. Lastly, the critical points computed by~\eqref{eq:critical_ystar} are $x^* = y^*=[\pm 0.236, \, 1.67]^\top$. Note that $I+\lambda^* L_s = \bsmtx 0 & 0 \\ 0 & -3 \esmtx$ is rank one. Therefore, the set of critical points forms a one-dimensional hypersphere, i.e., two endpoints of an interval in $\Rn{2}$.

\begin{figure}
    \centering
    \includegraphics[width=1\linewidth]{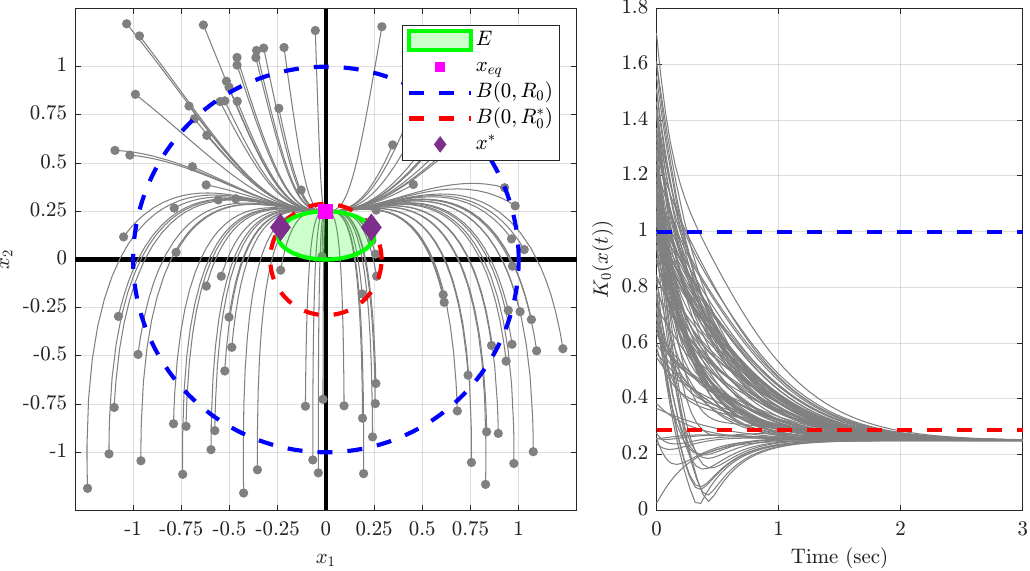}
    \caption{Illustrations of trapping regions of the bounded two-state system~\ref{eq:TwoStateSys_example}. The left plot shows the ellipsoid with non-decreasing energy $E$, valid trapping regions $B(0,R_0)$ and $B(0,R_0^*)$, and critical points $x^*$. Note that the radius $R_0$ (in blue) is an over-estimation discussed in Section~\ref{sect:prel_LorenzCond}, and $R_0^*$ (in red) is the tightest radius computed by the SDP-based method proposed in Section~\ref{sect:Verif_Computing_beta}. The critical points $x^*$ (purple diamonds) are computed by the process proposed in Section~\ref{sect:Verif_ystar}. The right plot shows the energy versus time evaluated along the state trajectories shown in the left plot. All energy trajectories decrease below $R_0^*$ in this system and well below $R_0$.
    }
    \label{fig:TwoStateSys_Example}
\end{figure}

Figure~\ref{fig:TwoStateSys_Example} illustrates these results. The left plot shows state trajectories (light gray) in the state space starting from random initial conditions. All trajectories eventually converge into both valid trapping regions $B(0, R_0)$ and $B(0, R_0^*)$. Note that the estimated radius $R_0$ is an over-estimation, and the radius $R_0^*$ is the tightest radius that contains the ellipsoid $E$ with nondecreasing energy~\eqref{eq:EnergyNonDecreasingSet} in green. Two critical points $x^*$ are marked by purple diamonds. The critical points have norm $R_0^*$ and correspond to the points in $E$ that are furthest from the origin. The right plot shows the energy versus time evaluated along the state trajectories shown in the left plot. All energy trajectories eventually decrease below $R_0^*$ and well below $R_0$, as shown by the analyses. Overall, the proposed SDP-based analysis provides a more accurate characterization of the trapping region than Theorem~\ref{theorem:TrappingRegion}.


Note that this system has a single equilibrium point at $x_{eq}:=(0, 0.25)$, and this is marked by a magenta rectangle on the left plot of Figure~\ref{fig:TwoStateSys_Example}. The energy $\frac{1}{2}\norm{x-x_0}^2$ is a Lyapunov function that proves this equilibrium point is globally asymptotically stable. This further implies that the system~\eqref{eq:TwoStateSys} is bounded. However, this is not true in general for lossless quadratic systems, i.e., they can have unstable equilibrium points and still have bounded trajectories (see next example). Moreover, the proposed SDP-based method can be used to verify the boundedness of lossless quadratic systems without directly computing any equilibrium point (as is done in Lyapunov analysis).

\subsection{Lorenz Chaotic Attractor}   \label{sect:example_Lorenz}

Consider the Lorenz system~\cite{lorenz1963deterministic}:
\begin{subequations} \label{eq:Lorenz}
    \begin{align}
        \frac{d x_1}{dt} &= -\sigma x_1 + \sigma x_2 \\
        \frac{d x_2}{dt} &= \rho x_1 - x_2 - x_1 x_3 \\
        \frac{d x_3}{dt} &= -\alpha x_3 + x_1x_2, 
    \end{align}
\end{subequations}
where $\sigma = 10, \rho = 28$ and $\alpha = \frac{8}{3}$. The system has no stable equilibrium points and exhibits chaotic behavior. Nonetheless, the system is shown to be bounded in Schlegel and Noack\cite{schlegel2015long}. Specifically, it can be verified analytically that $A_s(m)$ is negative definite with the translation $m= [ 0,\, 0, \, \rho + \sigma]^\top$. Moreover, a bound on the trapping region, given by Theorem~\ref{theorem:TrappingRegion}, is $B(m, R_m)$ with $R_m=101.33$. Note that the estimated radius $R_m$ is different from the radius reported in Schlegel and Noack~\cite{schlegel2015long}, as there is a minor typo in the estimation formula as pointed out in Kaptanoglu et al~\cite{Kaptanoglu2021}.

Next, we demonstrate the SDP-based results introduced in Section~\ref{sect:Verif}. First, the SDP~\eqref{SDP:min_lambda} was solved yielding $a^*=-1$ and $m=[0, \, 0, \, 38]^\top$. It follows from Theorem~\ref{theorem:TrappingRegion_SDP} and $a^*<0$ that a trapping region exists. Note that the solution $m$ is the same as the translation $m$ used in Schlegel and Noack~\cite{schlegel2015long}. This is not a typical case, as the set of $m$ such that $A_s(m)$ is negative definite is not unique, and the solution depends on the exact formulation of both approaches. Next, we solve the SDP~\eqref{eq:QCQPdual} and compute the smallest trapping region $B(m,R_m^*)$ with $R_m^*=39.25 < R_m = 101.33$. Lastly the critical points are characterized by~\eqref{eq:critical_ystar} as $x^*=y^*+m$ with $y^* = [0, \, \pm 24.82, \, -30.4]^\top$. 

\begin{figure}
\centering
\includegraphics[width=1\textwidth]{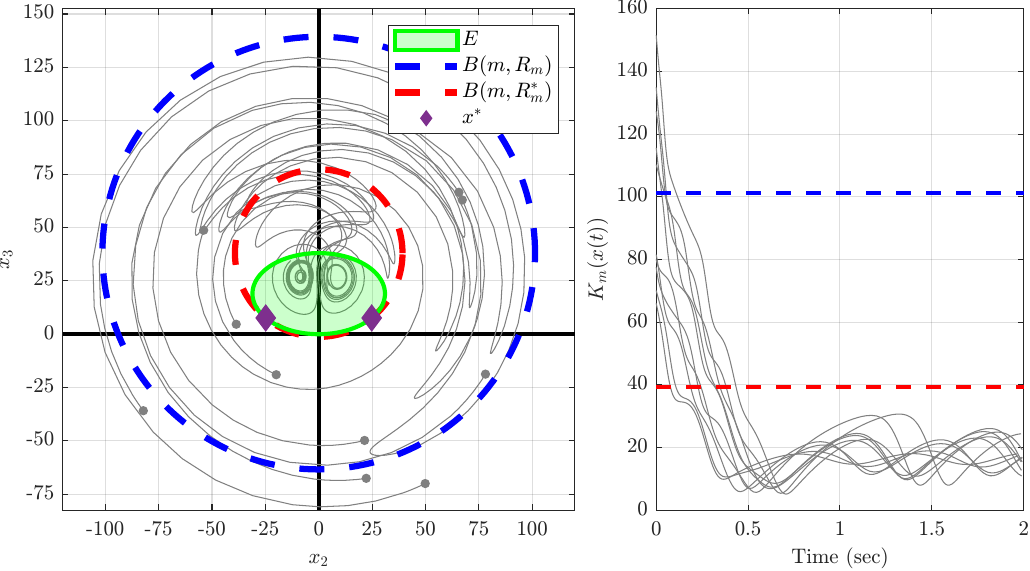}
\caption{Trapping regions of the Lorenz attractor computed by Theorem~\ref{theorem:TrappingRegion} (blue) and the proposed SDP-based analysis (red). The left plot visualizes trajectories in the $x_2-x_3$ plane with $x_1=0$. The green ellipsoid is the region where energy is instantaneously non-decreasing. The right plot shows the energy versus time evaluated along the state trajectories shown in the left plot. All trajectories are ultimately bounded in the region computed by the proposed SDP-based analysis (red), which captures the boundedness behavior more accurately than Theorem~\ref{theorem:TrappingRegion} (blue).
}
\label{fig:Lorenze}
\end{figure}

Figure~\ref{fig:Lorenze} illustrates these results. The left plot shows ten trajectories with randomly sampled initial conditions in the $x_2-x_3$ plane with $x_1=0$. Trajectories do not converge to any steady state or even a limit cycle since the system is chaotic. However, all trajectories are bounded as they are trapped in the boundedness region shown above. The red and blue circles are the computed boundedness region by our method and Theorem~\ref{theorem:TrappingRegion} respectively. The green ellipsoid is the area of instantaneously non-decreasing energy as characterized by the quadratic inequality in~\eqref{eq:EnergyNonDecreasingSet}. The purple diamonds mark the critical points $x^*$, which are the furthest points with non-decreasing energy from $m$. All trajectories converge to and stay inside of both blue and red circles as dictated by Theorem~\ref{theorem:TrappingRegion} and~\ref{theorem:TrappingRegion_SDP} respectively. The right plot presents the energy $K_m(x(t))$ evolution of the ten trajectories. As predicted, all energies eventually decay below and stay below levels of $R_m$ (blue) and $R_m^*$ (red). Furthermore, our analysis (red) provides a more accurate description of the boundedness behavior of the Lorenz system. 

\subsection{High-Dimensional Systems}


The example illustrates the proposed method on a higher dimensional system. The example is constructed by stacking multiple Lorenz systems and coupling them via random coordinate rotation. Specifically, consider a collection of $K$ systems with Lorenz dynamics~\eqref{eq:Lorenz}:
\begin{align*}
    \dot{z}^{(i)}(t) = L^{(i)} z^{(i)}(t) + f^{(i)} \left( z^{(i)}(t) \right) \,\,  i=1,\ldots K,
\end{align*}
where
\begin{align}
  L^{(i)} = \bmtx -\sigma & \sigma & 0 \\ \rho & -1 & 0 \\ 0 & 0 & -\alpha \emtx 
   \,\,  \mbox{ and } \,\,
  f^{(i)}( z^{(i)} ) = \bmtx 0 \\ - z^{(i)}_1 z^{(i)}_3 \\ z^{(i)}_1 z^{(i)}_2 \emtx.
\end{align}
We use the same parameters for each subsystem: $\sigma = 10$, $\rho =28$, and $\alpha = \frac{8}{3}$. 
Next, create a larger system by stacking these Lorenz subsystems: 
\begin{align*}
    \dot{z}(t) = L_z z(t)+ f_z( z(t) )
    \mbox{ where }
    z(t):= \bmtx z^{(1)}(t) \\ \vdots \\ z^{(K)}(t) \emtx \in \Rn{3K}.
\end{align*}
The matrix $L_z\in\Rn{3K\times 3K}$ and function $f_z:\Rn{3K} \to \Rn{3K}$ are obtained by appropriate concatenation of the Lorenz subsystems.  Here $f_z$ is a lossless quadratic function as each Lorenz subsystem has a lossless quadratic term. Finally, define the state for our high-dimensional example as $x=W z$ where $W\in \Rn{3K \times 3K}$ is an orthonormal matrix. The dynamics of this stacked and rotated system are:
\begin{align}   \label{eq:HighDimSys}
  \dot{x}(t) = L x(t) + f( x(t) )
  \,\mbox{ where } \,
   L = W L_z W^\top
   \, \mbox{ and }  \, 
   f(x) = W f_z( W^\top x). 
\end{align}   
The quadratic function in~\eqref{eq:HighDimSys} remains lossless after the rotation because:
\begin{align*}
    x^\top f(x)     =   z^\top f( z) = 0 
\end{align*}
The resulting high-dimensional system~\eqref{eq:HighDimSys} has a trapping region and is bounded by construction since each subsystem is bounded. However, the system~\eqref{eq:HighDimSys} is not obviously bounded as the dynamics of each Lorenz subsystem are coupled via the orthogonal transformation $x=Wz$.

The boundedness of~\eqref{eq:HighDimSys} is studied with 18 values of $K$ chosen on a logarithmic grid from $K=1$ to $K=223$. The corresponding state dimension of~\eqref{eq:HighDimSys} is $n=3K$ and this goes from 3 to 669. Ten systems with random orthogonal matrices $W$ are generated for each value of $K$. The SDP-based trapping region analysis in Theorem~\ref{theorem:TrappingRegion_SDP} is applied to each of these ten $3K$-dimensional systems~\eqref{eq:HighDimSys}. The analyses verified that the system~\eqref{eq:HighDimSys} has a trapping region $B(m, R_m^*)$ for each $K$ and each random $W$. The solution of the SDP~\eqref{SDP:min_lambda} yields $a^*=-1$ for each trial. This is the same result as for a single Lorenz system (Section~\ref{sect:example_Lorenz}) because the combined system~\eqref{eq:HighDimSys} is simply multiple copies of Lorenz dynamics in rotated coordinates. However, the coordinate shift $m$ is not unique since the first state in Lorenz subsystems has no nonlinearity. Any translation in this direction does not affect the shifted linear symmetric matrix $A_s(m)$. Hence, optimal solutions $m$ form a set and only one of the solutions is returned by the SDP solver. Furthermore, the radius $R_m^*$ of the combined system~\eqref{eq:HighDimSys} is not the same for each trial as the radius $R_m^*$ depends on the coordinate shift $m$. Nonetheless, all systems are bounded as expected.

\begin{figure}
    \centering
    \includegraphics[width=1\linewidth]{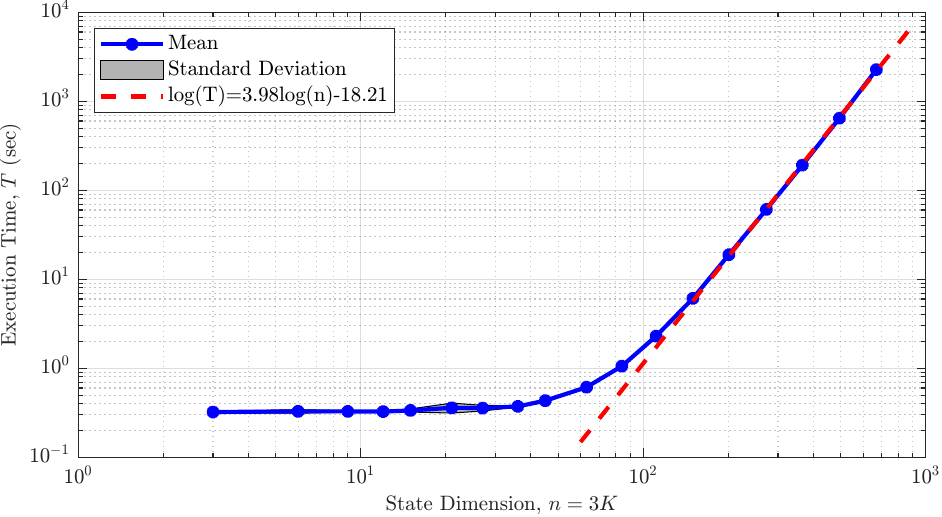}
    \caption{Computation times $T$ for the proposed SDP-based analysis versus the number of states $n=3K$ in the high-dimensional system~\eqref{eq:HighDimSys}. The computation includes solving the SDP~\eqref{SDP:min_lambda} and SDP~\eqref{eq:QCQPdual}. The mean (blue) increases as number of states grows, and the standard deviation (gray) is negligible for large state dimensions. For large state dimensions, the execution time scales with $O(n^{3.98})$ and is close to the theoretical complexity $O(n^4)$.
    }
    \label{fig:HighDimSys}
\end{figure}

Figure~\ref{fig:HighDimSys} illustrates the solve time with respect to the number of states $n=3K$. The computation time is the time to solve the two SDPs~\eqref{SDP:min_lambda} and~\eqref{eq:QCQPdual}. The results show that the computation time increases as $n$ increases (blue line). For a large number of states, the computation time scales empirically with $O(n^{3.98})$. This empirical estimate (red dashed line) was computed by linear regression on the data with $n>150$. This matches the theoretical complexity of the primal-dual interior-point SDP solver in MOSEK (Section 13.3 of MOSEK~\cite{mosek}). Specifically, the complexity of a primal-dual interior-point SDP solver is $\max{\{c_v c_d^3, \, c_v^2 c_d^2, \, c_v^3\}}$ (Chapter 11.8.3 of Boyd and Vandenberghe~\cite{boyd2004convex}), where $c_v$ is the number of decision variables and $c_d$ is the dimension of LMI. The first SDP~\eqref{SDP:min_lambda} has $(c_v,c_d) = (n+1, n)$ and the complexity is $O(n^4)$. The second SDP~\eqref{eq:QCQPdual} has $(c_v,c_d) = (2,n+1)$ and the complexity is $O(n^3)$. Hence, the overall theoretical complexity is $O(n^4)$ and is close to the observed complexity $O(n^{3.98})$. Note that the computation is dominated by the first SDP~\eqref{SDP:min_lambda} to verify the existence of a trapping region. Hence, the radius of the trapping region is relatively cheap to compute by the second SDP~\eqref{eq:QCQPdual} once the coordinate shift $m$ is found.  Note that the implementation is not optimized and the choice of $K=223$ is based on the memory size of our desktop platform. The analysis can scale to even larger systems by more efficient implementation, customized SDP solver, and improved computational hardware. 

\section{Conclusion} \label{sect:Conclusion}

This paper considers systems with lossless quadratic nonlinearities and proposes a convex optimization approach to analyze boundedness using trapping regions. The contribution is threefold. 
First, the necessary and sufficient condition for a trapping region to exist is formulated as a convex SDP. This SDP provides theoretical and numerical advantages to verify the condition. Readily available SDP solvers can verify or falsify the existence of a trapping region with guarantees. 
Second, the least conservative radius of the trapping region at a given coordinate, if one exists, is computed by solving a QCQP exactly using duality theory. The solution of QCQP also characterizes the critical states on the boundary of the trapping region with instantaneously non-decreasing energy. This approach reduces the conservative estimation compared to the literature and provides important testing scenarios for further attention. 
Lastly, numerical examples are presented to illustrate the analysis process and its computational performance. The proposed method can easily verify that a trapping region exists and provides the precise estimate of the region for systems up to $\sim O(100)$ states. 

The proposed convex-optimization-based method provides multiple additional opportunities for future work. 
The first direction is to incorporate this convex condition into modeling procedures. Several attempts~\cite{Kaptanoglu2021, goyal2023guaranteed} use the Schlegel and Noack trapping region condition (Theorem~\ref{theorem:TrappingRegion}) to propose system identification algorithms with a physics-informed prior. Our convex SDP (Theorem~\ref{theorem:TrappingRegion_SDP}) could potentially serve as a subroutine and improve the computational efficiency and modeling performance of such algorithms. 
The second direction is to use the convex condition in controller syntheses. SDPs are widely used to compute stabilizing controllers in both linear systems and nonlinear systems~\cite{boyd1994linear}. Notably, the convexity of static output feedback synthesis for lossless quadratic systems has been studied~\cite{MUSHTAQ2024ConvexSOF}. The trapping region condition being convex implies similar convexity structures potentially exist and can benefit controller synthesis tasks. 
The last direction here is the robustness analysis of the trapping region condition through convex optimization. Specifically, one could generalize from systems with lossless nonlinearities, as considered in this paper, to weakly lossless nonlinearities~\cite{Peng2023LocalStability}. Moreover, we could extend our convex condition to analyze the robustness of trapping regions by incorporating existing robust stability results~\cite{dullerud2013RobustControlCVX}. 
In summary, this convex formulation of the trapping region condition establishes potential connections to a large amount of literature in control theory utilizing convexity.

\bmsection*{Author contributions}
The authors confirm their contribution to the paper as follows: 
study conception and design: Shih-Chi Liao, A. Leonid Heide, Maziar S. Hemati, Peter J. Seiler; 
data collection: Shih-Chi Liao; 
analysis and interpretation of results: Shih-Chi Liao, Maziar S. Hemati, Peter J. Seiler; 
visualization: Shih-Chi Liao, A. Leonid Heide; 
draft manuscript preparation: Shih-Chi Liao, ;
manuscript revision: Maziar S. Hemati, Peter J. Seiler.
All authors reviewed the results and approved the final version of the manuscript. 

\bmsection*{Acknowledgments}
This material is based upon work supported by the Army Research Office under grant number W911NF-20-1-0156 and the Air Force Office of Scientific Research under grant number FA9550-21-1-0434. 
The authors acknowledge valuable discussions with Diganta Bhattacharjee. 

\bmsection*{Financial disclosure}
None reported.

\bmsection*{Conflict of interest}
The authors declare no potential conflict of interest.

\bmsection*{Supporting information}
Complete source code for the implementation and examples can be found on GitHub at \url{https://github.com/SCLiao47/Boundedness_LosslessQuadSys}.

\bmsection*{ORCID}
Shih-Chi Liao \orcidSCL{0000-0001-8311-3753} \\
A. Leonid Heide \orcidSCL{0009-0006-1522-245X} \\
Maziar S. Hemati \orcidSCL{0000-0003-1831-051X} \\
Peter J. Seiler \orcidSCL{0000-0003-3423-1109} \\




\bibliography{Bibliography}

\end{document}